\documentclass{amsart}

\usepackage{amsmath,amssymb,amsfonts,enumerate,amsthm, amscd}

\newcommand{\pd}{\mbox{pd}\,}
\newcommand{\id}{\mbox{id}\,}
\newcommand{\fd}{\mbox{fd}\,}
\newcommand{\gd}{\mbox{Gdim}\,}

\newtheorem{theorem}{Theorem}[section]
\newtheorem{lemma}[theorem]{Lemma}
\newtheorem{proposition}[theorem]{Proposition}
\newtheorem{corollary}[theorem]{Corollary}

\theoremstyle{definition}

\theoremstyle{remark}
\newtheorem{remark}[theorem]{Remark}
\newtheorem{example}[theorem]{Example}

\theoremstyle{Definition and Notation}

\begin{document}
\bibliographystyle{amsplain}

\title[]{Gorenstein global dimension of Semi-primary rings}

\author{Mohammed Tamekkante}
\address{Mohammed Tamekkante\\Department of Mathematics, Faculty of Science and Technology of Fez, Box 2202, University S.M. Ben Abdellah Fez, Morocco.}



\begin{abstract}
The aim of this paper is the study of Gorenstein global and weak
dimensions of semi-primary rings.
\end{abstract}
\maketitle
\section{Introduction}
Throughout the paper, all rings are associative with identity, and
all modules are unitary.\\
Let $R$ be a ring, and let $M$ be an $R$-module. As usual we use
$\pd_R(M)$, $\id_R(M)$ and $\fd_R(M)$ to denote, respectively, the
classical projective dimension, injective dimension and flat
dimension of $M$.\\

For a two-sided Noetherian ring $R$, Auslander and Bridger
\cite{Aus bri} introduced the $G$-dimension, $\gd_R (M)$, for
every finitely generated $R$-module $M$. They showed that there is
an inequality $\gd_R (M)\leq \pd_R (M)$ for all finite $R$-modules
$M$, and equality holds if $\pd_R (M)$ is finite.

Several decades later, Enochs and Jenda \cite{Enochs,Enochs2}
defined the notion of Gorenstein projective dimension
($G$-projective dimension for short), as an extension of
$G$-dimension to modules that are not necessarily finitely
generated, and the Gorenstein injective dimension ($G$-injective
dimension for short) as a dual notion of Gorenstein projective
dimension. Then, to complete the analogy with the classical
homological dimension, Enochs, Jenda and Torrecillas \cite{Eno
Jenda Torrecillas} introduced the Gorenstein flat dimension. Some
references are
 \cite{Bennis and Mahdou2, Christensen, Christensen
and Frankild, Enochs, Enochs2, Eno Jenda Torrecillas, Holm}.\\

Recall that a left (resp., right) $R$-module $M$ is called
Gorenstein projective if, there exists an exact sequence of
projective left (resp., right) $R$-modules:
$$\mathbf{P}:...\rightarrow P_1\rightarrow P_0\rightarrow
P^0\rightarrow P^1\rightarrow ...$$ such that $M\cong
Im(P_0\rightarrow P^0)$ and such that the operator $Hom_R(-,Q)$
leaves $\mathbf{P}$ exact whenever $Q$ is a left (resp., right)
projective $R$-module. The resolution $\mathbf{P}$ is called a complete projective resolution. \\
The left and right Gorenstein injective $R$-modules are defined
dually.\\
 And an $R$-module $M$ is called left (resp., right)
Gorenstein flat if, there exists an exact sequence of flat left
(resp., right) $R$-modules:
$$\mathbf{F}:...\rightarrow F_1\rightarrow F_0\rightarrow
F^0\rightarrow F^1\rightarrow ...$$ such that $M\cong
Im(P_0\rightarrow P^0)$ and such that the operator $I\otimes_R-$
(resp., $-\otimes_RI$) leaves $F$ exact whenever $I$ is a right
(resp., left) injective $R$-module. The resolution $\mathbf{F}$ is called complete flat resolution.\\
The  Gorenstein projective, injective and flat
dimensions are defined in term of resolution and denoted by $Gpd(-)$, $Gid(-)$ and $Gfd(-)$ respectively (please see \cite{Christensen, Enocks and janda, Holm}).\\

In the rest of this papers, the word $R$-module will mean left
$R$- module unless explicitly stated otherwise.\\
In \cite{Bennis and Mahdou2}, the authors prove the equality
$$sup\{Gpd_R(M)|M\;is\;an\;R-module\}=sup\{Gid_R(M)|M\;is\;an\;R-module\}$$
They called the common value of the above quantities the left
Gorenstein global dimension of $R$ and denoted it by
$l.Ggldim(R)$. Similarly, they set
$$l.wGgldim(R)=\{Gfd_R(M)|M\;is\;an\;R-module\}$$
which they called the left Gorenstein weak dimension of $R$.

 Recall
that a ring $R$ is called semi-primary if there is a two-sided
nilpotent ideal $N$ of $R$, which we call the radical of $R$, such
that $R/N$ is semi-simple. It is clear that if $R$ is
semi-primary, its radical is unique (see \cite{Nicolson}). Some
familiars examples of semi-primary rings are:

\begin{example}
In the following cases the rings $R$ are semi-primary with radical
$J$:
\begin{enumerate}
    \item $R=\left(%
\begin{array}{cc}
  F & F \\
  0 & F \\
\end{array}%
\right)$ where $F$ is a field,  $J=\left(%
\begin{array}{cc}
  0 & F \\
  0 & 0 \\
\end{array}%
\right)$ with $J^2=0$.
 \item $R=K[X]/(X^2)$ where $K$ is field, $J=(\overline{X})$ and
 $J^2=0$.
  \item Camillo Example (\cite[Example 2.6]{Nicolson}).
\end{enumerate}
\end{example}
For more different examples (non-Noetherian, non-commutative, non
self-injective...) of semi-primary rings please see
(\cite{Nicolson}).

 The  purpose of this
papers is to characterize the
 left  Gorenstein global and weak dimensions  of  semi-primary rings.

\section{Main results}

We begin with the first main result which give a characterization
of the left Gorenstein
global dimension of  semi-primary rings.\\
Recall that the word $R$-module will mean left $R$- module unless
explicitly stated otherwise.

\begin{theorem}\label{theorem1}

Let $R$ be a semi-primary ring with radical $N$. Then,\\

 $\begin{array}{cccc}
  a) & lG.gldim(R) & = & Gpd_R(R/N) \\
  b) &  & = & Gid_R(R/N) \\
   c) &  & = & sup_C\;Gpd_R(C) \\
   d) & &=  & sup_C\;Gid_R(C) \\
\end{array}$\\

where $C$ ranges ranges over all simple left $R$-modules.\\
If
therefore $R$ is not  quasi-Frobenius ring then,\\

$\begin{array}{cccc}
   e) & lG.gldim(R) & = & 1+ Gpd_R(N) \\
 \end{array}$
\end{theorem}
To prove this theorem, we need the following Lemma.
\begin{lemma}\label{lemma1}
Let $0\rightarrow N\rightarrow N' \rightarrow N'' \rightarrow 0$
be an exact sequence of $R$-modules. Then,
\begin{enumerate}
  \item $Gpd_R(N')\leq max\{Gpd_R(N),Gpd_R(N")\}$ with equality
  if $Gpd_R(N'')\neq Gpd_R(N)+1$.
  \item $Gpd_R(N'')\leq max\{Gpd_R(N'),Gpd_R(N)+1\}$ with equality
  if $Gpd_R(N')\neq Gpd_R(N)$.
\end{enumerate}
\end{lemma}
\begin{proof}
Using \cite[Theorems 2.20 and 2.24]{Holm} the argument is
analogous to the one of \cite[Corollary 2, p. 135]{bourbaki}.
\end{proof}
\begin{proof}[Proof of Theorem \ref{theorem1}]

$(a).$ Obviously, by definition, we have $lG.gldim(R)\geq
Gpd_R(R/N)$. Then, only the other inequality need a proof and we
may assume $Gpd_R(R/N)<\infty$. We claim $Gpd_R(M)\leq Gpd_R(R/N)$
for each  $R$-module $M$.  Firstly, let $M$ be a  $R$-module such
that $NM=0$. Note that $M$ can be considered a left module over
the semi-simple ring $R/N$ by setting $\bar{x}.m=x.m$, for each
$m\in M$ and $x\in R$. Clearly this modulation is well defined
since
$$\bar{x}=\bar{y}\quad \Longrightarrow \quad x-y\in N\quad
\Longrightarrow \quad (x-y).m=0$$
  Then, $M$ is a left projective
$R/N$-module (recall that  $R/N$ is semi-simple). Therefore, it is
a direct summand of  a left free $R/N$-module $(R/N)^{(I)}$. Since
$NM=0$, we can consider $M$ as a direct summand of $(R/N)^{(I)}$
as an $R$-modules. Then, from
\cite[Propostoion 2.19]{Holm}, $Gpd_R(M)\leq Gpd_R((R/N)^{(I)})=Gpd_R(R/N)$. \\
Now, let $M$ be an arbitrary  $R$-module and let $k$ be the
smaller positive integer such that $N^kM=0$ (such integer  exists
since $N$ is nilpotent). Consider the family of short exact
sequences of  $R$-modules: $$0\rightarrow N^{k-i+1}M\rightarrow
N^{k-i}M \rightarrow N^{k-i}M/N^{k-i+1}M\rightarrow 0$$ where
$0<i\leq k$. Then, by Lemma \ref{lemma1}(1) we have
$$Gpd_R(N^{k-i}M)\leq
sup\{Gpd_R(N^{k-i+1}M),Gpd_R(N^{k-i}M/N^{k-i+1}M)\}$$ But,
$Gpd_R(N^{k-i}M/N^{k-i+1}M)\leq Gpd_R(R/N)$ since
$N(N^{k-i}M/N^{k-i+1}M)=0$. Thus, $Gpd_R(N^{k-i}M)\leq
sup\{Gpd_R(N^{k-i+1}M),Gpd_R(R/N)\}$. So, we conclude that
$Gpd_R(M)\leq sup\{Gpd_R(N^{k-1}M,Gpd_R(R/N)\}$. Again we have
$Gpd_R(N^{k-1}M)\leq Gpd_R(R/N)$ since $N(N^{k-1}M)=0$. Hence,
$Gpd_R(M)\leq Gpd_R(R/N)$, as desired.\\

$(b).$ Similarly to $(a)$ it suffices  to prove that for every
left  $R$-module $M$ such that $NM=0$ we have $Gid_R(M)\leq
Gid_R(R/N)$. The rest of the proof is the same lines as $(a)$ by
using the dual of Lemma \ref{lemma1}. Let $M$ be such module.
Then, $M$ is a direct summand of a left free $R/N$-module
$(R/N)^{(I)}$. If we identify $M$ to a submodule of $(R/N)^{(I)}$
we get
$$M\subseteq (R/N)^{(I)}\subseteq \Pi_{I}(R/N)$$ Then $M$ is a
direct summand of $\Pi_{I}(R/N)$ (as an $R/N$-modules and also as
an $R$-modules) since $M$ is an injective $R/N$-module (since
$R/N$ is semi-simple). Then, Using the injective version of
\cite[Proposition 2.19]{Holm}, we have $$Gid_R(M)\leq
Gid_R(\Pi_{I}(R/N)=Gid_R(R/N)$$ as desired.\\

$(c).$ Since $R/N$ is semi-simple, $R/N\cong \oplus C_i$, finite
direct sum of simple left $R$-modules, where the $C_i$ have the
property that if $C$ is a left simple $R$-module, then $C\cong
C_i$ for some $i$. Therefore, by \cite[Proposition 2.19]{Holm},
$sup_C\{Gpd_R(C)\}=lGpd_R(R/N)=lG.gldim(R)$, where $C$ ranges
over all left simple $R$-modules.\\

$(d).$ Since the direct sum $R/N\cong \oplus C_i$ is finite, we
can replace the direct sum by the direct product and use the
injective version of \cite[Proposition 2.19]{Holm}. Thus, $(d)$ is
proved in an analogous fashion to $(c)$. \\

$(e).$ Suppose that  $R$ is not  quasi-Frobenius. Then
$lG.gldim(R)=Gpd_R(R/N)>0$ (\cite[Proposition 2.6]{Bennis and
Mahdou2}). Therefore, from Lemma \ref{lemma1}(2) we  deduce from
the exact
 sequence $$0\longrightarrow N \longrightarrow R \longrightarrow R/N \longrightarrow 0$$
that $Gpd_R(R/N)=1+ Gpd_R(N)$.
\end{proof}

The next Proposition give a functorial description of the the left
Gorenstein global dimension of semi-primary rings provided this
value  is finite.
\begin{proposition}
Let $R$ a semi-primary rings with radical $N$ and with finite left
Gorenstein global dimension and let $n>0$ be an integer. The
following are equivalents:
\begin{enumerate}
    \item $lG.gldim(R)<n$,
    \item $Ext_R^n(I,R/N)=0$ for every injective $R$-module $I$,
    \item $Ext_R^n(I,C)=0$ for every injective $R$-module $I$ and every simple left $R$-module $C$,
    \item $Tor_R^n(R/N,I)=0$ for every injective $R$-module $I$ (and $R/N$ is consider as a right R-module),
    \item $Tor_R^n(C,I)=0$ for every injective $R$-module $I$ and every simple right $R$-module $C$,
    \item $Ext_R^n(R/N,P)=0$ for every projective $R$-module $P$,
    and
    \item $Ext_R^n(C,P)=0$ for every projective $R$-module $I$ and every simple left $R$-module
    $C$.
\end{enumerate}
\end{proposition}
\begin{proof}
 By \cite[Theorem 1.1]{Bennis and
Mahdou2}, we have $$lG.gldim(R)=sup\{Gid_RM|M\; an\; R-module\}.$$
Then, by \cite[Theorem 2.22]{Holm}, $lG.gldim(R)<n$ if, and only
if,\ $Ext^i(I,M)=0$ for each $i\geq n$ and every  $R$-module $M$
and every  injective $R$-module $I$. Hence we conclude that
$$lG.gldim(R)<n\quad \Longleftrightarrow\quad pd_R(I)<n\; for\;
every\; injective \; R-module\; I.$$ So, using
 \cite[Proposition 7]{Auslander}, we have the equivalence of
 $(1)$, $(2)$, $(3)$, $(4)$ and $(5)$.\\
 Using \cite[Lemma 2.1]{Bennis and Mahdou2} and \cite[Proposition
 10]{Auslander} we obtain the equivalence of $(1)$, $(6)$ and
 $(7)$.
\end{proof}

Now, we give our second main result which characterize the left
Gorenstein weak dimension of coherent semi-primary rings:
\begin{theorem}\label{theorem2}
Let $R$ be a right coherent semi-primary ring with radical $N$.
Then,
$$l.wGgldim(R)=Gfd_R(R/N)=sup_C\;Gfd_R(C)$$

where $C$ ranges ranges over all simple left $R$-modules.
\end{theorem}

To prove this Theorem we need the following Lemma:
\begin{lemma}\label{lemma2}
Let $0\rightarrow N\rightarrow N' \rightarrow N'' \rightarrow 0$
be an exact sequence of modules over a right coherent ring $R$.
Then, $Gfd_R(N')\leq max\{Gfd_R(N),Gfd_R(N)\}$ with equality
  if $Gfd_R(N'')\neq Gfd_R(N)+1$.
\end{lemma}
\begin{proof}
Using \cite[Theorem 3.15]{Holm} and \cite[Theorem 3.14]{Holm} the
proof is similar to \cite[Corollary 2, p. 135]{bourbaki}.
\end{proof}
\begin{proof}[Proof of Theorem \ref{theorem2}]
Using \cite[Proposition 3.13]{Holm} and Lemma \ref{lemma2},
 the proof is the same lines as of proof of the equality $(a)$ and $(c)$ of
 Theorem \ref{theorem1}.
 \end{proof}

 The next Proposition is an application  of Theorem
 \ref{theorem1}:
 \begin{proposition}\label{proposition}
 Let $R$  be a semi-primary ring such that each simple left
$R$-module is isomorphic to a left ideal in $R$, then:
$\{gldim(R),l.Ggldim(R)\}\in \{0,\infty\}$.
\end{proposition}
\begin{proof} The classical result $gldim(R)\in\{0,\infty\}$ is
exactly \cite[Proposition 14]{Auslander} and the Gorenstein
version is proved by the same way. For exactness we give the proof
here. Suppose that $lG.lgidm(R)=n$, $0<n<\infty$. By Theorem
\ref{theorem1} we have $lG.gldim(R)= Gpd_R(C)$, where $C$ is a
simple left $R$-module. By hypothesis, $C \cong I$, where $I$ is
an ideal in $R$. Thus $Gpd_R(I)=n$. Consider the exact sequence of
$R$-modules:
$$0\rightarrow I\rightarrow R \rightarrow R/I\rightarrow 0$$ Since
$n> 0$, $R/I$ is not Gorenstein projective (\cite[Theorem
2.5]{Holm}). Therefore by Lemma \ref{lemma1} $Gpd_R(R/I)=1+n$.
Contradiction with the fact that $Gpd_R(R/I)\leq l.Gglim(R)=n$.
This cotradiction finish the proof.
\end{proof}
\begin{remark}[Proposition 15, \cite{Auslander}]\label{remark} The hypothesis of Proposition \ref{proposition} is satisfied in each of the
following cases:
\begin{enumerate}
    \item $R$ is a direct sum of a finite number of primary rings (a semi-primary
ring $R$ is primary if $R/N$ is a simple ring).
    \item $R$ is a semi-primary commutative ring.
    \item $R$ is a quasi-Frobenius ring (i.e; Noetherian and self-injective ring).
\end{enumerate}
\end{remark}
\begin{corollary}
Every commutative semi-primary rings with finite Gorenstein global
dimension is quasi-Frobenius.
\end{corollary}
\begin{proof} This Corollary is a direct consequence of
Proposition \ref{proposition}, Remark \ref{remark} and
\cite[Proposition 2.6]{Bennis and Mahdou2}.
\end{proof}

\bibliographystyle{amsplain}

\begin{thebibliography}{10}


\bibitem{Auslander}M. Auslander, \textit{On the dimension of modules and algebras (III)}, global
dimension, Nagoya Math. J., 9 (1955), 67-77.
\bibitem{Aus bri}M. Auslander and M. Bridger;  \textit{Stable module theory},
Memoirs. Amer. Math. Soc., 94, American Mathematical Society,
Providence, R.I., 1969.


\bibitem{Bennis and Mahdou2}D. Bennis and N. Mahdou; \textit{Global
Gorenstein Dimensions}, accepted for publication in Proceedings of
the American Mathematical Society. Available from
math.AC/0611358v4 30 Jun 2009.


\bibitem{bourbaki}   N. Bourbaki, \textit{Alg$\grave{e}$bre Homologique}, Chapitre 10, Masson, Paris,
(1980).


\bibitem{Christensen} L. W. Christensen; \textit{Gorenstein dimensions}, Lecture Notes in Math., Vol. 1747, Springer,
Berlin, (2000).
\bibitem{Christensen and Frankild}
L. W. Christensen, A. Frankild, and H. Holm; \textit{On Gorenstein
projective, injective and flat dimensions - a functorial
description with applications}, J. Algebra 302 (2006), 231-279.
\bibitem{Enochs} E. Enochs and O. Jenda; \textit{On Gorenstein injective modules}, Comm.
Algebra 21 (1993), no. 10, 3489-3501.

\bibitem{Enochs2} E. Enochs and O. Jenda; \textit{Gorenstein injective and projective
modules}, Math. Z. 220 (1995), no. 4, 611-633.
\bibitem{Eno Jenda Torrecillas}
E. Enochs, O. Jenda and B. Torrecillas; \textit{Gorenstein flat
modules}, Nanjing Daxue Xuebao Shuxue Bannian Kan 10 (1) (1993)
1-9.
\bibitem{Enocks and janda} E. E. Enochs and O. M. G. Jenda, \textit{Relative Homological Algebra}, de Gruyter Expositions in
Mathematics, Walter de Gruyter and Co., Berlin, 2000.


\bibitem{Holm} H. Holm; \textit{Gorenstein homological dimensions}, J. Pure Appl.
Algebra 189 (2004), 167-193.



\bibitem{Nicolson} W. K. Nicholson and M. F. Youssif; \textit{Quasi-Frobenius Rings}, Cambridge University
Press, vol. 158, 2003.







\end{thebibliography}

\end{document}